\documentclass[letterpaper, 10 pt, conference]{ieeeconf}  

\IEEEoverridecommandlockouts                              
\usepackage{graphics} 
\usepackage{epsfig} 
\usepackage{mathptmx} 
\usepackage{times} 
\usepackage{amsmath} 
\usepackage{amssymb}  
\usepackage{subcaption}
\usepackage{cite}
\let\labelindent\relax
\usepackage{enumitem}

\title{\LARGE \bf A Branch-Decomposition Approach to Power Network Design}

\author{Kin Cheong Sou}

\begin{document}

\maketitle
\thispagestyle{empty}
\pagestyle{empty}

\begin{abstract}
This paper proposes a procedure to solve combinatorial power network design problems such as phasor measurement unit (PMU) placement and protection assignment against cyber-physical attacks. The proposed approach tackles the design problems through solving a dominating set problem on the power network graph. A combined branch-decomposition and dynamic programming procedure is applied to solve the dominating set problem. Contrary to standard integer programming and heuristic/evolutionary approaches for power network design, the proposed approach is exact and requires only polynomial computation time if the graph is planar and has small branch-width. A planarization technique is explored for problem instances with nonplanar graphs. A case study in this paper with benchmark power networks verifies that planarity and small branch-width are not uncommon in practice. The case study also indicates that the proposed method shows promise in computation efficiency.
\end{abstract}

\section{Introduction}
Our society depends heavily on the proper operation of network systems including intelligent transport systems, electric power distribution and transmission systems etc. These systems are supervised and controlled through Supervisory Control And Data Acquisition (SCADA) systems. For instance, in the electric power transmission grid, SCADA systems collect measurements through remote terminal units (RTUs) and send them to the state estimator to estimate the system states. The estimated states are used for subsequent operations such as contingency analysis (for system health monitoring) and optimal power flow dispatch (for control). Any malfunctioning of these operations can lead to significant social and economical consequences such as the northeast US blackout of 2003 \cite{US_blackout2003}.

Because of its importance, the SCADA measurement system has been the subject of extensive studies. For electric power network, observability analysis has traditionally been an important topic in measurement system analysis (e.g., \cite{wu1985network, krumpholz1980power, castillo2006observability, Abur_Exposito_SEbook,Monticelli_SEbook}). Observability analysis determines if the state of the system can be deduced from measurements. Another important measurement system related topic which has attracted a lot of recent attention is cyber-physical security. One of the purposes of cyber-physical security is to analyze various types of data attacks and their consequences on the system (e.g.,\cite{LRN09,dan2010stealth,sandberg2010security,bobba2010detecting,kosut2010malicious,giani2011smart,kim2011strategic}).

Both observability and cyber-physical security have their design aspects, on which this paper focuses. A design question for observability seeks the most economical setup of measurement system to ensure observability (e.g., \cite{gou2000direct}). In some cases, the system must be robust in the sense that even if one or more measurements are missing the system should remain observable \cite{xu2004observability}. More recent work on observability-enabling design focuses on the strategic placement of phasor measurement units (PMUs), which are becoming more and more prevalent in modern power networks (e.g., \cite{xu2004observability, aminifar2010contingency, chakrabarti2009placement, 6563646, Anderson201451}). A version of the PMU placement problem, known as the power dominating set problem, is extensively studied in computer science related communities. A wealth of theoretical results are available (e.g., \cite{haynes2002domination,dorbec2013generalized,guo2005improved,brueni2005pmu}). On the other hand, for cyber-physical security the design typically seeks the minimum cost protection strategy so that, according to the chosen attack and defense model, no data attack is possible (e.g., \cite{dan2010stealth, bobba2010detecting, kim2011strategic}).

Because of the combinatorial nature of the PMU placement problems, solution procedures roughly fall into two groups: integer programming (e.g., \cite{xu2004observability, aminifar2010contingency, chakrabarti2009placement}) and heuristic evolutionary algorithms (e.g., \cite{peng2006optimal}) as surveyed in \cite{manousakis2012taxonomy}. A similar categorization can be made for the methods to design protection against cyber-physical attacks. The available design methods enjoy obvious advantages such as ease to implement and formulate, modeling flexibility and in the case of integer programming, guaranteed optimality (if given enough time) and highly optimized solvers such as CPLEX and Gurobi. However, there are also shortcomings associated with the available methods. For instance, a fundamental limitation for integer programming is that its computation requirement (time and memory) increases sharply with problem instance size. For evolutionary algorithms, neither optimality nor computation time can be guaranteed. One of the main purposes of this paper is to investigate an alternative design methodology with a different set of advantages and disadvantages. The proposed method is based on a graph decomposition technique called branch-decomposition \cite{ROBERTSON1991153}. In fact, branch-decomposition has been applied with success to solve difficult problems (e.g., \cite{christian2002linear, cook2003tour, bian2008algorithms,arnborg1991easy,bodlaender1994tourist}). A branch-decomposition generates a hierarchy of expanding subgraphs to enable streamlined dynamic programming to solve difficult (e.g., NP-hard) combinatorial graph optimization problems vital to the design problems considered in this paper (i.e., PMU placement and cyber-physical protection). With suitable problem structure such as planarity and small branch-width (to be described in sequel), the branch-decomposition based approach can solve the design problems in polynomial time with no ``enormous hidden constant''. This overcomes the limitation of integer programming and heuristic methods. Indeed, the possibility to exploit problem structure is also a distinguishing factor between the branch-decomposition approach and the available, ``general-purpose'' approaches. In addition, branch-decomposition can be ``re-used''. Once a branch-decomposition is constructed, it can be re-used to efficiently solve many different combinatorial problems for the same graph. This is analogous to LU factorization in equation solving with different right-hand-sides.
 
At its core, the proposed design methodology in this paper sets up and solves a (minimum) dominating set problem on the given power network graph. Given a graph, a dominating set is a subset of the vertex set such that every vertex not in the dominating set is adjacent to at least one vertex in the dominating set. A minimum cardinality dominating set is sought in the dominating set problem. There exist a large collection of results on the dominating set problem (e.g., \cite{haynes1998fundamentals,haynes1998domination}). The results most related to this paper are those on planar dominating set problem (e.g., \cite{fomin2006dominating, alber2002fixed, demaine2005fixed}). However, these results are of theoretical nature and focus on establishing general properties which could be too conservative for practical purposes. Instead, this paper focuses on applying some of these algorithms in practical settings to make observations which are difficult to obtain in a general theoretical setting. The same theoretical versus practical perspectives distinguishes this paper from the results on power dominating set problem \cite{haynes2002domination,dorbec2013generalized,guo2005improved,brueni2005pmu}. Of possible exception is \cite{guo2005improved}, which applies tree-decomposition and dynamic programming to solve the power dominating set problem. However, practically branch-decomposition and tree-decomposition can be different because there is ``no enormous hidden constant'' with branch-width computation \cite{fomin2006dominating}.

{\bf Outline:} in Section~\ref{sec:applications}, the considered network design problems will be described, and their relationship to the dominating set problem will be explained. In Section~\ref{sec:computing_MDS} selected methods for computing minimum dominating set will be described. Section~\ref{sec:BD} describes the first step of the proposed branch-decomposition based method. That is, branch-width and branch-decomposition and their computation will be described. Section~\ref{sec:DP} explains the second step of the proposed method. It explains how dynamic programming can be applied to find a minimum dominating set, when a branch-decomposition is given. Section~\ref{sec:algo_summary} provides a short summary of the entire proposed methodology. Section~\ref{sec:case_study} describes the results of a numerical case study on computing minimum dominating sets for well-known power network benchmarks. The case study compares the accuracy and time-efficiency between the proposed method, the integer programming approach and a heuristic greedy algorithm.

\section{Application motivations} \label{sec:applications}
A power network can be modeled as a graph where the vertices are buses, and the edges are transmission lines. Typically, in order to ensure proper operation of the power network, the voltage phasors at the buses need to be measured directly or determined indirectly through measurements. The optimization for the setup and protection of the measurement system often gives rise to the dominating set problem. 

\subsection{PMU placement problem} \label{sec:PMU}
In a version of the PMU placement problem in \cite{xu2004observability}, the goal is to measure the voltage phasors at all buses using PMUs. In the problem setup, if a PMU is located at a bus, the voltage phasor at the bus can be measured directly. In addition, the voltage phasors at all buses adjacent to the bus with PMU can be measured indirectly. Hence, to measure all voltage phasors in the network, it is necessary and sufficient that the set of buses with PMUs form a dominating set in the network. Assuming that the total cost of the PMU-based measurement system is proportional to the number of PMUs installed, the PMU placement problem can then be written as a dominating set problem. A PMU is installed at a bus if and only if it is included in the minimum dominating set.

Another version of PMU placement problem is called the power dominating set problem (e.g., \cite{haynes2002domination}). Similar to the first version, a PMU at a bus can measure the voltage phasors at the bus and all its adjacent buses. In addition, it is assumed that the power injections at all buses are known (e.g., load, generator, zero-injection buses). Therefore, if the voltage phasors at a bus and all but one adjacent buses are measured (directly or indirectly), then the phasor at the remaining bus is indirectly measured (i.e., Kirchhoff's current law) \cite{PDS_IP,Kneis2006145}. A power dominating set is a subset of buses with PMU installed so that all voltage phasors are measured according to the above rules. The problem in \cite{haynes2002domination} seeks the minimum cardinality power dominating set. By the measurement rule, a (standard) dominating set of the network is a power dominating set, even though the converse is not true. As a result, a solution to the dominating set problem provides an approximate solution to the minimum power dominating set problem.

\subsection{Perfect protection problem} \label{sec:protection}
The measurement model in \cite{dan2010stealth} is the DC power flow model \cite{Abur_Exposito_SEbook,Monticelli_SEbook}. In the setup of \cite{dan2010stealth}, each bus is equipped with a remote terminal unit (RTU) which can measure all or part of the following quantities: (net) active power injection at the bus and active power flows on the transmission lines incident to the bus. However, the RTU cannot measure the voltage phasor directly as in the case of PMU. Instead, the voltage phasors are estimated based on the measurements from the RTUs.

In power system operations, it is a standard practice to employ a ``bad data detection'' scheme to detect possible anomalies in the measurements. Typically if there is only one attacked measurement, the detection scheme can detect it. However, coordinated additive attacks on multiple measurements locations may evade detection. Additivity means the attacked measurement is the sum of the original measurement and a modification (the attack). Reference \cite{LRN09} reports that the attacks will avoid detection if they can be interpreted as voltage-phasor-induced active power flows and injections. That is, detection is avoided if some fictitious voltage phasors can be associated with the buses such that (a) the modification to the active power flow on a transmission line is proportional (with appropriate constant) to the difference of the fictitious voltage phasors at the terminal buses, and (b) at each bus the modification to active power injection satisfies Kirchhoff's current law with the modifications to active power flows on incident transmission lines. To counter these attacks, \cite{dan2010stealth} considers the scenario in which buses can be protected by installing encryption devices at RTUs. Once a bus is protected, all measurements related to the RTU on that bus is considered attack-free. Reference \cite{dan2010stealth} poses the problem of perfect protection as follows: protect the minimum number of buses with the constraint that no detection-avoiding attack can be staged. If there is a bus, denoted $i$, and all its adjacent buses are unprotected, then the active power injection measurement at $i$ and active power flow measurements on all transmission lines incident to $i$ can be attacked. These attack modifications, according to (a) and (b), can indeed be made undetected since they can be induced by the following fictitious voltage phasors: zero at all buses except at bus $i$. Hence, as argued in \cite{dan2010stealth}, a necessary condition for perfect protection is that the protected buses must form a dominating set of the network. While a dominating set need not imply perfect protection, its computation is necessary since it is 
the initial guess for the perfect protection algorithm in \cite{dan2010stealth}.

\subsection{Computing minimum dominating set} \label{sec:computing_MDS}
In general, the dominating set problem is NP-hard \cite{garey2002computers} and is not fixed parameter tractable \cite{downey1999parameterized}. In this paper, three different solution approaches for the dominating set problem are considered: 
\begin{enumerate}
\item Integer programming exact solution
\item Max-degree greedy approximate algorithm
\item Branch-decomposition based (exact or approximate) solution
\end{enumerate}
The integer programming approach finds a minimum dominating set as the solution of the following integer program:
\begin{equation} \label{opt:DS_IP}
\begin{array}{cl}
\underset{x}{\text{minimize}} & \sum\limits_{i = 1}^{|V|} x_i \vspace{1mm} \\
\text{subject to} & G + I_{|V|} \ge {\bf 1} \vspace{1mm} \\
& x_i \in \{0,1\},
\end{array}
\end{equation}
where $G \in \{0,1\}^{|V| \times |V|}$ is the adjacency matrix of graph $(V,E)$, and ${\bf 1}$ denotes the vector of ones. The greedy algorithm studied in this paper is as follows:

\vspace{3mm}
\noindent Greedy algorithm for approximate minimum dominating set \\
\noindent Input: Graph $(V,E)$ \\
Output: Dominating set $D$
\begin{itemize}
\item Initialize $D \leftarrow \emptyset$
\item While $D$ is not dominating, do
\begin{itemize}
\item $v \leftarrow$ vertex in $V \setminus D$ with the maximum number of neighbors either not in $D$ or not adjacent to $D$ (i.e., not dominated by $D$)
\item $D \leftarrow D \cup \{v\}$
\end{itemize}
End (of While)
\end{itemize}
The third approach, the branch-decomposition based approach, can be divided into two steps to be discussed in the next two sections:

\vspace{3mm}
\noindent {\bf Algorithm 1:} Branch-decomposition based method for planar dominating set problem
\begin{description}
\item[Step 1 (Section~\ref{sec:BD}):] Compute the branch-width using the ratcatcher algorithm (e.g., \cite{seymour1994call}), and construct an optimal branch-decomposition (e.g., \cite{gu2008optimal}).
\item[Step 2 (Section~\ref{sec:DP}):] Solve the dominating set problem using dynamic programming \cite{fomin2006dominating}, whose computation is organized by the optimal branch-decomposition.
\end{description}



\section{Branch-decomposition} \label{sec:BD}
The first step of the proposed branch-decomposition based method for dominating set problem is to compute the branch-width and an optimal branch decomposition of the input graph. Two notations are needed: for any graph $G$, $V(G)$ denotes the vertex set and $E(G)$ denotes the edge set. Next, we define branch-decomposition and branch-width as follows: Given a graph $(V,E)$, a \emph{branch-decomposition} \cite{ROBERTSON1991153} is a pair $(T, \tau)$ where $T$ is a unrooted binary tree with $|E|$ leaf vertices, and $\tau$ is a bijection from $E$ to the set of leaf vertices of $T$. Every non-leaf vertex of $T$ has degree three. For any edge $e \in E(T)$, the subgraph $(V(T), E(T) \setminus \{e\})$ has two connected components denoted $T_1(e)$ and $T_2(e)$. Let $E_1(e)$ and $E_2(e)$ be the edge subsets of $E$ corresponding to the leaf vertices of $T_1(e)$ and $T_2(e)$, respectively. Let $\omega(e) := \{v \in V \; \vline \;  v \in e_1, v \in e_2, \; \text{for some $e_1 \in E_1(e)$ and $e_2 \in E_2(e)$} \}$. That is, $\omega(e)$ is the subset of $V$ whose members are incident to edges in both $E_1(e)$ and $E_2(e)$. Then, the width of branch-decomposition $(T, \tau)$ is $\max_{e \in E(T)} | \omega(e) |$. The \emph{branch-width} of graph $(V,E)$ is the minimum width over all branch-decompositions of $(V,E)$. A branch-decomposition of graph $(V,E)$ is optimal if its width is the branch-width of the graph.

\begin{figure}[!tbh]
\begin{center}
 \includegraphics[width=70mm]{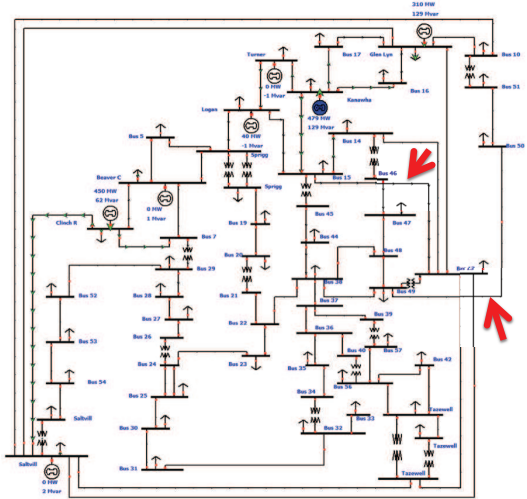}
\end{center}
 \caption{A drawing of the IEEE 57-bus benchmark power network \cite{IEEE57_CSL}. From this drawing, it can be seen that removing only two edges $\{13,15\}$, $\{49,50\}$ (with red arrows pointed) will render the remaining subgraph planar.}
 \label{fig:IEEE57}
\end{figure}

In general, it is NP-hard to compute the branch-width of a graph \cite{seymour1994call}. However, if restricted to \emph{planar} graphs computing branch-width becomes polynomial time \cite{seymour1994call} (see, for example, \cite{diestel2005graph} for the definition of planar graph). For instance, the ratcatcher algorithm in \cite{seymour1994call} computes the branch-width of a planar graph $(V,E)$ in $O(\log_2(|V|) (|V|+|E|)^2)$ time. Given the branch-width, constructing an optimal branch-decomposition for a planar graph is also polynomial time (e.g., $O(|V|^3)$ in \cite{gu2008optimal}). It is one of the purposes of this paper to report that a significant number of power network graphs (e.g., \cite{IEEE_power_grid_benchmarks}) are indeed planar and ``almost planar'' as in Fig.~\ref{fig:IEEE57}. Nevertheless, not all power network graphs are planar. For example, the IEEE 118-bus graph is planar but the IEEE 300-bus graph is not. To address the planarity issue, this paper suggests that the input graph is first tested for planarity (e.g., Boyer-Myrvold algorithm \cite{diestel2005graph}). If the graph is planar, proceed with the second step dynamic programming in Algorithm 1 (cf.~Section~\ref{sec:computing_MDS}, Section~\ref{sec:DP}). If the graph is not planar, then compute an edge-maximal planar subgraph as follows:

\vspace{3mm}
\noindent Algorithm to compute edge-maximal planar subgraph \\
\noindent Input: Connected graph $(V,E)$ \\
Output: Edge-maximal planar subgraph $(V,E_p)$
\begin{enumerate}
\item Find a spanning tree $T$ of $(V,E)$; initialize $E_p \leftarrow E(T)$.
\item Test if there exists $e \in E \setminus E_p$ such that $(V, E_p \cup \{e\})$ is planar.
\item If there exists $e$ as in 2), update $E_p \leftarrow E_p \cup \{e\}$; goto 2). If no $e$ exists, goto 4).
\item Return $(V,E_p)$ as the desired edge-maximal planar subgraph of $(V,E)$.
\end{enumerate}
Notice that we consider edge-maximal planar subgraph instead of edge-maximum planar subgraph, because finding the later is NP-hard \cite{garey2002computers}. Once an edge-maximal planar subgraph is found, proceed with the second step in Algorithm 1 with this subgraph. The minimum dominating set on the subgraph is a suboptimal dominating set on the original graph. For the test cases examined in this paper, this planarization-based technique is found to be sufficient (see Section~\ref{sec:case_study}). However, more sophisticated planarization techniques are available (e.g., \cite{gutwenger2010application}). Also note that for real-world power network examples, there are inherent drawings (crossings allowed) of the underlying graphs corresponding to the physical locations of the buses and transmission lines. These drawings can be useful in identifying good edge-maximal planar subgraphs. See Fig.~\ref{fig:IEEE57} for an example of the IEEE 57-bus benchmark.

\section{Dynamic programming based on branch-decomposition} \label{sec:DP}
A branch-decomposition generates a hierarchy of expanding subsets of edges, and correspondingly an expanding set of edge-induced subgraphs. This hierarchy enables the second step of the proposed methodology because it constitutes an optimal substructure for dynamic programming to solve the dominating set problem. The description in this section is based on \cite{fomin2006dominating} and \cite{Gu880}. We make the following definitions:
\begin{description}
\item[$T'$:] Let $(T, \tau)$ be an optimal branch-decomposition of graph $(V,E)$. Construct $T'$, a rooted tree from $T$ by inserting two new vertices: (a) vertex $z$ into any edge $\{u,v\} \in E(T)$, and (b) vertex $r$ (the root vertex) forming an edge $\{z,r\}$. Specifically, let $\{u,v\} \in E(T)$, then $T' = \big(V(T) \cup \{z,r\}, (E(T) \setminus \{u,v\}) \cup \{\{u,z\}, \{z,v\}, \{z,r\}\}\big)$. The rooted tree $T'$ is used to organize the computations in dynamic programming.
\item[$G_e$:] For each $e \in E(T')$, a leaf vertex $u \in V(T')$ is a descendant of $e$ if the (unique) path from root $r$ to $u$ traverses $e$. Let $V_{T'}(e)$ denote the subset of leaf vertices of $T'$ that are descendants of $e$. Then, $G_e$ is defined as the edge-induced subgraph of $(V,E)$ whose edge set is $\tau^{-1}(V_{T'}(e))$. Note that if $e, e' \in E(T')$ are such that the path from $r$ to $e$ traverses $e'$ then $G_e \subseteq G_{e'}$. Traversing $T'$ from root to leaves leads to a hierarchy $\{G_e \; \vline \; e \in E(T')\}$ of subgraphs, in which a member $G_{e'}$ includes another member $G_e$ if $e$ is ``below'' $e'$ in $T'$. This hierarchy leads to an optimal substructure for dynamic programming.
\item[$\omega'$:] Recall, in the definition of branch-width in Section~\ref{sec:BD}, the definition of $\omega(e) \subseteq V$ for $e \in E(T)$. We define a function $\omega' : E(T') \mapsto V$ as follows: let $\omega'(\{z,r\}) = \emptyset$, $\omega'(\{u,z\}) = \omega'(\{z,v\}) = \omega(\{u,v\})$ and $\omega'(e) = \omega(e)$ for all $e \in E(T') \setminus \{\{u,z\}, \{z,v\}, \{z,r\}\}$. Notice that $\omega'(e) \subseteq V(G_e)$ because all members of $\omega'(e)$ are incident to $\tau^{-1} (V_{T'}(e)) = E(G_e)$, and the subset inclusion can be strict. The vertex subset $\omega'(e)$ can be interpreted as the ``boundary'' of $V(G_e)$, and $\omega'(e)$ defines the ``state-space'' for dynamic programming, for ``stage $e$''.
\end{description}
For each $e \in E(T')$, we color the vertices of $\omega'(e) \subseteq V$ in three possible colors:
\begin{description}
\item[Black:] represented by 1, meaning that the vertex is in the dominating set.
\item[White:] represented by 0, meaning that the vertex is not in the dominating set but it is dominated by at least one vertex in $V(G_e)$ that is in the dominating set.
\item[Grey:] represented by $\hat{0}$, meaning that the vertex is not in the dominating set. It is not necessarily dominated inside $G_e$.
\end{description}
The symbol $c \in \{0,\hat{0},1\}^{|\omega'(e)|}$ denotes a specific coloring of all vertices in $\omega'(e)$. $c$ encodes the vector of a particular ``state'' in dynamic programming. Associated with $c$, we define
\begin{description}
\item[$A_e(c)$:] For each $e \in E(T') \setminus \{z,r\}$, we define function $A_e : \{0, \hat{0}, 1\}^{|\omega'(e)|} \mapsto \mathbb{N} \cup \infty$ such that for each coloring $c \in \{0,\hat{0},1\}^{|\omega'(e)|}$, $A_e(c)$ is defined as $|D_e|$, where $D_e \subseteq V(G_e)$ has minimum cardinality among all subsets of $V(G_e)$ satisfying
\begin{itemize}
\item All black vertices in $\omega'(e)$ are in $D_e$.
\item No white vertex in $\omega'(e)$ is in $D_e$, but a white vertex must be adjacent to a vertex in $D_e$.
\item No grey vertex in $\omega'(e)$ is in $D_e$.
\item Each vertex in $V(G_e) \setminus \omega'(e)$ is either in $D_e$ or adjacent to a vertex in $D_e$.
\end{itemize}
If, for a particular coloring $c \in \{0,\hat{0},1\}^{|\omega'(e)|}$, there is no subset of $V(G_e)$ satisfying the above four conditions then $A_e(c) := \infty$.
\end{description}
Note that $A_e$ plays the role of ``cost-to-go'' function, at ``stage $e$'', in the standard dynamic programming terminology. $\{0,\hat{0},1\}^{|\omega'(e)|}$ is the set of all possible states at stage $e$ and the cardinality of the state-space is $3^{|\omega'(e)|} \le 3^{\text{BW}}$. Therefore, for dynamic programming to be efficient, the branch-width of graph $(V,E)$ should be small. The edges $e \in E(T')$ incident to a leaf of $T'$ can be thought of as the ``final stages''. All values of $A_e$ for a leaf-incident edge $e$ can be enumerated since there are at most $3^2 = 9$ choices (corresponding to the coloring of two end vertices of a single edge $\tau^{-1}(V_{T'}(e)) \in E$). For all edges $e \in E(T')$ not incident to a leaf in $T'$, $A_e$ is assembled (through a dynamic programming recursion) from the associated functions $A_{e_1}$ and $A_{e_2}$, where $e_1$ and $e_2$ are two children edges of $e$ defined as follows:
\begin{description}
\item[$e_1$, $e_2$:] For each edge $e \in E(T')$ not incident to a leaf vertex in $T'$, let $v_e \in V(T')$ denote the vertex incident to $e$ and having $e$ on the path from the root $r$ to $v_e$. There are two other edges incident to $v_e$. Let $e_1$ and $e_2$ denote these two edges. They are the two children of $e$.
\end{description}
To assemble $A_e$ from $A_{e_1}$ and $A_{e_2}$, more definitions are needed:
\begin{description}
\item[$X_1, X_2, X_3, X_4$:] \begin{displaymath}
\begin{array}{l}
X_1 := \omega'(e) \setminus \omega'(e_2), \vspace{0.5mm} \\
X_2 := \omega'(e) \setminus \omega'(e_1), \vspace{0.5mm} \\
X_3 := \omega'(e) \cap \omega'(e_1) \cap \omega'(e_2), \vspace{0.5mm} \\
X_4 := (\omega'(e_1) \cup \omega'(e_2)) \setminus \omega'(e).
\end{array}
\end{displaymath}
\end{description}
For each coloring $c$ for vertices in $\omega'(e)$, only some colorings $c_1$ and $c_2$ for $\omega'(e_1)$ and $\omega'(e_2)$ (respectively) are consistent with $c$ and profitable to consider in dynamic programming recursion:
\begin{itemize}
\item For $u \in X_1$, $c(u) = c_1(u)$. In this case, for each $v \in V$ adjacent to $u$, $v \in V(G_e)$ if and only if $v \in V(G_{e_1})$.
\item For $u \in X_2$, $c(u) = c_2(u)$. In this case, for each $v \in V$ adjacent to $u$, $v \in V(G_e)$ if and only if $v \in V(G_{e_2})$.
\item For $u \in X_3$, if $c(u) \in \{\hat{0},1\}$ then $c(u) = c_1(u) = c_2(u)$. If $c(u) = 0$ then either (a) $c_1(u) = 0, c_2(u) = \hat{0}$ or (b) $c_1(u) = \hat{0}, c_2(u) = 0$. In fact $c(u) = \hat{0}$ is consistent with $(c_1(u), c_2(u)) \in \{0,\hat{0}\}^2$. However, only $c_1(u) = c_2(u) = \hat{0}$ is most profitable. If $c(u) = 0$ then $u$ is dominated by at least one $v \in V(G_e)$ adjacent to $u$. $v$ must be in $G_{e_1}$, $G_{e_2}$ or both. Hence, at least one of $c_1(u)$ and $c_2(u)$ must be $0$. Again, the choice of $c_1(u) = c_2(u) = 0$ is not profitable and hence it is excluded.
\item For $u \in X_4$, exactly one of the three cases must hold: (a) $c_1(u) = c_2(u) = 1$, (b) $c_1(u) = 0, c_2(u) = \hat{0}$ or (c) $c_1(u) = \hat{0}, c_2(u) = 0$. Vertices in $X_4$ are internal vertices of $V(G_e)$. They do not have adjacent vertices outside $V(G_e)$. Therefore, if $u$ is not in the dominating set there must be a vertex in $V(G_{e_1})$ or $V(G_{e_2})$ (as they partition $V(G_e)$) in the dominating set. The choice $c_1(u) = c_2(u) = \hat{0}$ can therefore be ``too optimistic'' and it is excluded. On the other hand, the choice $c_1(u) = c_2(u) = 0$ is excluded since it is not the most profitable.
\end{itemize}
Then, for $e \in E(T')$ not incident to a leaf vertex, the dynamic programming recursion to assemble $A_e(c)$ from $A_{e_1}$ and $A_{e_2}$ is
\begin{equation} \label{eqn:Jmin_general}
\begin{array}{ccl}
A_e(c) & = & \min \Big\{A_{e_1} (c_1) + A_{e_2} (c_2) - \#_1(X_3, c_1) - \#_1(X_4, c_1) \vspace{0.5mm} \\ & & \vline \; \text{$c_1$, $c_2$ form $c$} \Big\},
\end{array}
\end{equation}
where $\#_1(X_3, c_1)$ denotes the number of vertices in $X_3$ colored 1 by coloring $c_1$ (note that for $c_1$, $c_2$ forming $c$, it holds that $\#_1(X_3, c_1) = \#_1(X_3, c_2)$). Similarly, $\#_1(X_4, c_1)$ denotes the number of vertices in $X_4$ colored 1 by coloring $c_1$ (and also by $c_2$). By definition of $\omega'$, there is no vertex in exactly one of $\omega'(e)$, $\omega'(e_1)$ and $\omega'(e_2)$. Hence, $X_3$ and $X_4$, which are disjoint, partition the set $\omega'(e_1) \cap \omega'(e_2) = V(G_{e_1}) \cap V(G_{e_2})$. Therefore, in (\ref{eqn:Jmin_general}) the two minus terms prevent double-counting of dominating set members in $\omega'(e_1) \cap \omega'(e_2)$.

For the root-incident edge $\{z,r\}$, $\omega'(\{z,r\}) = \emptyset$. Consequently, $X_1 = X_2 = X_3 = \emptyset$. Let $e_{r_1}$ and $e_{r_2}$ denote the two children edges of $\{z,r\}$. Then, $X_4 = \omega'(e_{r_1}) = \omega'(e_{r_2})$ since there is no vertex in exactly one of $\omega'(e_{r_1})$ or $\omega'(e_{r_2})$. The cardinality of the minimum dominating set of $(V,E)$ is
\begin{equation} \label{eqn:Jmin_root}
\underset{c_{r_1}, \; c_{r_2}}{\min} \;\; A_{e_{r_1}} (c_{r_1}) + A_{e_{r_2}} (c_{r_2}) - \#_1(\omega'(e_{r_1}), c_{r_1}),
\end{equation}
subject to the constraint that for all $u \in \omega'(e_{r_1})$, exactly one of the three conditions holds: (a) $c_{r_1}(u) = c_{r_2}(u) = 1$, (b) $c_{r_1}(u) = 0, c_{r_2}(u) = \hat{0}$ or (c) $c_{r_1}(u) = \hat{0}, c_{r_2}(u) = 0$.

Let $D^\star$ denote the minimum dominating set to be returned. Then the members of $D^\star$ can be decided by a tree traversal (e.g., depth-first) starting from the root: let $c^\star_{r_1}$, $c^\star_{r_2}$ be a minimizing pair in (\ref{eqn:Jmin_root}). Subsequently, for each $e \in E(T')$ such that $c^\star = \underset{c}{\text{argmin}} \; A_e(c)$ is known, we handle two cases:
\begin{itemize}
\item If $e$ is incident to a leaf vertex then for each $u \in \omega'(e)$, $c^\star(u) = 1$ implies $u \in D^\star$ and $u \notin D^\star$ if $c^\star(u) \in \{0, \hat{0}\}$. If $|\omega'(e)| = 1$ and suppose $u \in \omega'(e)$ and $v \notin \omega'(e)$, then $v \in D^\star$ if and only if $u \notin D^\star$.
\item If $e$ is not incident to a leaf vertex, then the optimal colorings $c^\star_1$ and $c^\star_2$ for the two children edges $e_1$ and $e_2$ can be read off from (\ref{eqn:Jmin_general}), with $c$ fixed as $c^\star$.
\end{itemize}
The process continues until all edges in $T'$ have been visited. 

Reference \cite{fomin2006dominating} shows that the time-complexity of the dynamic programming procedure is $O(3^{1.5 \text{BW}} |E(T')|) = O(3^{1.5 \text{BW}} |E|)$ because $|E(T')| = O(|E|)$. Therefore, it is crucial that the graph considered has \emph{small branch-width} in order for the dynamic programming to execute efficiently (with small branch-width the runtime is linear). One of the purposes of this paper is to experimentally demonstrate that a large number power network graphs \cite{IEEE_power_grid_benchmarks}, even if they are reasonably large, have very small branch-widths.

\section{Summary of the proposed procedure} \label{sec:algo_summary}
Summarizing the discussions in Section~\ref{sec:BD} and Section~\ref{sec:DP}, the steps of the proposed dominating set problem solution procedure is

\vspace{3mm}
\noindent {\bf Algorithm 1p:} Branch-decomposition based method for dominating set problem
\begin{description}
\item[Step 0:] Test if graph $(V,E)$ is planar (e.g., Boyer-Myrvold algorithm \cite{diestel2005graph}). If $(V,E)$ is planar, set $E_p = E$. If not planar, compute $(V,E_p)$ as an edge-maximal planar subgraph of $(V,E)$ with $E_p \subset E$.
\item[Step 1:] Compute the branch-width of planar graph $(V,E_p)$ using the ratcatcher algorithm (e.g., \cite{seymour1994call}), and construct an optimal branch-decomposition $(T,\tau)$ (e.g., \cite{gu2008optimal}).
\item[Step 2:] Use $(T,\tau)$ to set up the dynamic programming in Section~\ref{sec:DP} to solve the dominating set problem on $(V,E_p)$.
\item[Step 3:] Return $D_p$ as a dominating set of $(V,E)$. $D_p$ is optimal if $(V,E)$ is planar, but it is only suboptimal if $(V,E)$ is not planar.
\end{description}
Again, planarity (or ``almost planarity'') and small branch-width are crucial for the success of the proposed algorithm. The case study in the next section indicates that this is not uncommon in practice. 

\section{Numerical case study} \label{sec:case_study}

Graphs from the benchmark data set \cite{IEEE_power_grid_benchmarks,zimmerman2011matpower} are considered. The minimum dominating sets of the considered graphs are computed by solving integer program (\ref{opt:DS_IP}). To solve (\ref{opt:DS_IP}), two solvers are used -- IBM CPLEX (via MATLAB interface) and the branch-and-bound implementation in YALMIP called BNB \cite{YALMIP}. CPLEX is used because it is one of the standards in integer programming solvers, and BNB is used to examine the time-performance of a free but not-necessarily-optimized solver. In addition, the proposed branch-decomposition based procedure (Algorithm 1p in Section~\ref{sec:algo_summary}) is used to compute the minimum dominating sets (suboptimal in case the graph is not planar). To compute the branch-width and branch-decomposition, we utilize the state-of-the-art implementation by Prof.~Gu's group \cite{Bian2016156}. The dynamic programming computations, as described in Section~\ref{sec:DP}, are implemented in MATLAB (for ease of implementation). For certain graph operations (e.g., planarity test, finding a planar embedding), MatlabBGL is used \cite{Gleich06contentsmatlab}. In addition, for comparison, the greedy algorithm described in Section~\ref{sec:computing_MDS} is implemented and tested. The computations for branch-width and branch-decomposition are performed on a Linux machine with a 3 GHz CPU and 4GB of RAM. All other computations are performed on a Mac machine with a 2.5 GHz CPU and 8GB of RAM.

The test results of the studied graphs are shown in Table~\ref{tab:domination_number}. The meaning of each column in Table~\ref{tab:domination_number} is as follows: $|V|$ is the number of vertices. $|E|$ is the number of edges. $\text{BW}_p$ is the branch-width of the edge-maximal planar subgraph used to compute branch-decomposition for each test case (for a planar graph the edge-maximal planar subgraph is the graph itself). $|D^\star|$ is the domination number (i.e., cardinality of the minimum dominating set), as obtained by solving (\ref{opt:DS_IP}) using CPLEX. $|D^\star_{BD}|$ is the cardinality of the dominating set obtained using the proposed Algorithm 1p. In case of planar graph this cardinality is the same as the domination number, otherwise it is an upper bound of the domination number. $|D^\star_{gr}|$ is the cardinality of the dominating set returned by the greedy algorithm.
\begin{table}[h]
\begin{center}
\caption{Branch-width and domination number}
\begin{tabular}{|c|c|c|c|c|c|c|c|}
\hline
Name & \scriptsize{$|V|$} & \scriptsize{$|E|$} & planar & $\text{BW}_p$ &  \scriptsize{$|D^\star|$} &\scriptsize{$|D^\star_{BD}|$} & \scriptsize{$|D^\star_{gr}|$}  \\[2pt]
\hline
IEEE 9 & 9 & 9 & yes & 2 & 3 & 3 & 3\\
\hline
IEEE 14 & 14 & 20 & yes & 2 & 4 & 4 & 5\\
\hline
IEEE 24 & 24 & 34 & no & 3 & 7 & 7 & 9 \\
\hline
IEEE 30 & 30 & 41 & yes & 3 & 10 & 10 & 10 \\
\hline
IEEE 39 & 39 & 46 & yes & 3 & 13 & 13 & 14 \\
\hline
IEEE 57 & 57 & 78 & no & 4 & 17 & 17 & 21 \\
\hline
IEEE 118 & 118 & 179 & yes & 4 & 32 & 32 & 38 \\
\hline
IEEE 300 & 300 & 409 & no & 4 & 87 & 91 & 96\\
\hline
\end{tabular}
\label{tab:domination_number}
\end{center}
\end{table}
From Table~\ref{tab:domination_number} it can be seen that the branch-widths of the planar graphs used to construct the dominating sets are indeed small, meaning that Algorithm 1p is time-efficient for the cases listed. In addition, except for the IEEE 300 case Algorithm 1p returns dominating sets whose cardinalities are equal to the domination numbers. It is also seen that the proposed procedure is more accurate than the greedy algorithm.

The elapsed time for various parts in the case study is listed in Table~\ref{tab:time}. The units of all recordings in the table are in seconds. The meaning of each column in Table~\ref{tab:time} is as follows: BW+BD denotes the time to compute the branch-width and construct an optimal branch-decomposition of the graph. DP denotes the time to solve the dominating set problem using dynamic programming, when an optimal branch-decomposition is given. BW+BD+DP denotes the sum of the time for BW, BD and DP. This amounts to almost all computation time required by Algorithm 1p. The column BNB displays the time that the YALMIP branch-and-bound implementation uses to find minimum dominating sets. In the last row, $>3300$ means that the computation was aborted after 3300 seconds (when the 3000-iteration limit was reached with about 9.37\% optimality gap). Finally, the last column CPLEX denotes the time for CPLEX to solve the dominating set problem.
\begin{table}[h]
\begin{center}
\caption{Computation time (sec) for the tested algorithms}
\begin{tabular}{|c|c|c|c|c|c|}
\hline
Name & BW+BD & DP & \tiny{BW+BD+DP} & \scriptsize{BNB} & \scriptsize{CPLEX} \\
\hline
IEEE 9 & 0.03 & 0.066 & 0.096 & 0.31 & 0.0016 \\
\hline
IEEE 14 & 0.05 & 0.16 & 0.21 & 0.29 & 0.0066 \\
\hline
IEEE 24 & 0.13 & 0.51 & 0.64 & 0.39 & 0.0083 \\
\hline
IEEE 30 & 0.06 & 0.57 & 0.63 & 0.65 & 0.0063 \\
\hline
IEEE 39 & 0.09 & 0.54 & 0.63 & 0.31 & 0.0065 \\
\hline
IEEE 57 & 0.12 & 2.14 & 2.26 & 6 & 0.014 \\
\hline
IEEE 118 & 1 & 5.47 & 6.47 & 35 & 0.01 \\
\hline
IEEE 300 & 3.4 & 10.5 & 13.9 & $>3300$ & 0.024 \\
\hline
\end{tabular}
\label{tab:time}
\end{center}
\end{table}
Table~\ref{tab:time} indicates that the current implementation of Algorithm 1p is on par in time with BNB with smaller instances. For the largest instance studied (IEEE 300), Algorithm 1p is much more efficient than BNB mainly because BNB runtime scales badly with instance size. On the other hand, CPLEX, the industry standard, outperforms the current implementation of Algorithm 1. In order for the proposed branch-decomposition based method to be competitive, the implementation of dynamic programming should be more efficient (e.g., C instead of MATLAB implementation). The computation of the branch-decomposition can also be improved.

\begin{figure}[!tbh]
\begin{center}
 \includegraphics[width=70mm]{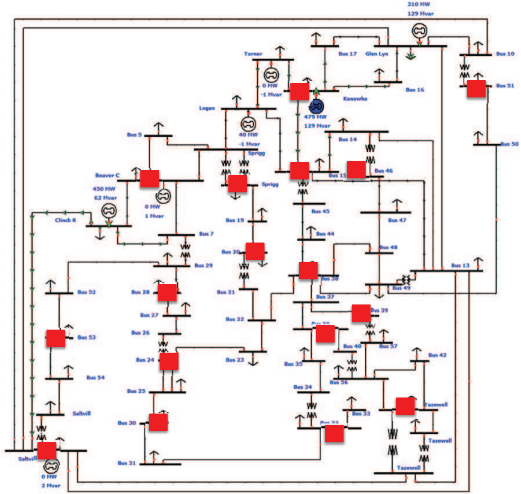}
\end{center}
 \caption{Optimal placement of PMUs for the IEEE 57-bus benchmark. Red rectangles indicate PMUs.}
 \label{fig:IEEE57_mds}
\end{figure}

Finally, to illustrate Algorithm 1p as a means to solve the PMU design problem, Fig.~\ref{fig:IEEE57_mds} shows the layout of an optimal placement of PMUs (red rectangles) for the IEEE 57-bus benchmark.

\section{Conclusion and discussion}
The case study in this paper indicates that graphs modeling practical power networks can possess nice properties such as small branch-width and planarity. In addition, even if the graphs are nonplanar, they tend to be ``almost planar'' as also remarked in \cite{brueni2005pmu}. This can be an inherent property (or tendency) due to economic considerations. These nice properties translate into the possibility to utilize branch-decomposition to efficiently solve many difficult combinatorial graph problems. This is showcased by the dominating set problem modeling the power network design in this paper. The dominating set problem, however, is only a start of the investigation into branch-decomposition enabled dynamic programming methodologies for other power network combinatorial problems. Currently, a not-so-naive implementation of the branch-decomposition base approach can be on par with a not-necessarily-optimized branch-and-bound implementation for solving the considered network design problems. In fact, the experimental result suggests that the branch-decomposition approach manages better with problem instance size, as it should. Nevertheless, contrasted with industry standards such as CPLEX, the branch-decomposition approach should be improved before it can become a serious competitor to the standard integer programming based approaches. Another direction for improvement would be to compute the branch-width and branch-decomposition even when the input graph is not planar. This problem is NP-hard, but practically possible method may exist for power network graphs in practice.

\section*{Acknowledgement}
The author is grateful to Dr.~Qianping Gu for the fruitful communications regarding branch-decomposition.

\bibliographystyle{IEEEtran}
\bibliography{bd_acc16}

\end{document}